\documentclass{article}

\usepackage{amsthm}
\usepackage{amsmath}
\usepackage{amssymb}
\usepackage{graphicx}
\usepackage{epsfig}

\newtheorem{theorem}{Theorem}

\newtheorem{definition}[theorem]{Definition}
\newtheorem{example}[theorem]{Example}
\newtheorem{remark}[theorem]{Remark}

\begin{document}

\begin{center} \LARGE{\textbf{
VARIETIES WITH DEGENERATE GAUSS MAPS WITH MULTIPLE FOCI AND
TWISTED CONES}}

\vspace*{3mm}

{\large  Maks A. Akivis and  Vladislav V. Goldberg}

\end{center}

Abstract. {\footnotesize The authors study in detail new types of
varieties with degenerate Gauss maps: varieties with multiple foci
and their particular case, the so-called twisted cones. They prove
an existence theorem for twisted cones and describe their
structure.

\vspace*{2mm}


\vspace*{2mm}


}

 \vspace*{3mm}

\setcounter{equation}{0}

\setcounter{section}{-1}

\section{Introduction}
A smooth $n$-dimensional variety $X$ of a projective space $P^N$
is called  {\em tangentially degenerate} or {\em a variety with a
degenerate Gauss map} if the rank of its Gauss mapping $\gamma: X
\rightarrow G (n, N)$ is less than $n, \; 0 \leq r = \mbox{{\rm
rank}} \; \gamma < n$. Here $x \in X, \; \gamma (x) = T_x (X)$,
and $T_x (X)$ is the tangent subspace to $X$ at $x$ considered as
an $n$-dimensional projective space $P^n$.  The number  $r$ is
also called the {\em rank} of $X, \; r = \mbox{{\rm rank}} \; X$.
The case $r = 0 $ is trivial one: it gives just an $n$-plane.

Let $X \subset  P^N$ be  an $n$-dimensional smooth variety with a
degenerate Gauss map. Suppose that  $0 < \mbox{{\rm rank}} \,
\gamma = r < n$. Denote by $L$ a leaf of this map, $L =
\gamma^{-1} (T_x) \subset X;\, \dim L = n - r = l$. The number $l$
is called the Gauss defect of the variety $X$ (see [FP 01], p. 89,
or [L 99], p. 52) or the index of relative nullity of $X$ (see [CK
52]).

A variety with a degenerate Gauss map of rank $r$ foliates into
their leaves  $L$ of dimension $l$, along which the tangent
subspace $T_x (X)$ is fixed. The foliation on $X$ with leaves $L$
is called the {\em Monge--Amp\`{e}re foliation} (see, for example,
[D 89] or [I 98, 99b]).

However, unlike a traditional definition of the foliation, the
leaves of the Monge--Amp\`{e}re foliation  have singularities.
This is a reason that in general its leaves are not diffeomorphic
to a standard leaf. We assume that singular points belong to
the leaf $L$, and hence the leaf is an $l$-dimensional subspace
of the space $P^N$.

The tangent subspace $T_x (X)$ is fixed when a point $x$ moves
along  regular points of $L$. This is the reason that we denote it
by $T_L, \, L \subset T_L$. A pair $(L, T_L)$ on $X$ depends on
$r$ parameters.

The varieties of rank $r < n$ are multidimensional analogues of
developable surfaces of a three-dimensional Euclidean space.
They  were first considered by \'{E}. Cartan [C 16] in connection
with his study of metric deformations of hypersurfaces, and in
[C 19] in connection with his study of manifolds of constant curvature.
Recently varieties with degenerate Gauss maps of rank $r < n$
are intensively studied both from the projective point of view and
the Euclidean point of view.

The main  results on the geometry of varieties with degenerate
Gauss maps and further references can be found in Chapter 4 of the
book [AG 93] and in the recently published paper [AG 01b].

Griffiths and Harris [GH 79] (Section  2, pp. 383--393) considered
varieties with degenerate Gauss maps from the point of view
of algebraic geometry. Following [GH 79], Landsberg  published notes
[L 99] which are in some sense an  update to the   paper [GH 79].
Section 5 (pp. 47--50) of these notes is devoted to varieties with
degenerate Gauss maps.

In particular, in [GH 79] Griffiths and Harris  presented a
structure theorem for  varieties with degenerate Gauss maps. They
asserted that such varieties are ``built up from cones and
developable varieties'' (see [GH 79], p. 392). They gave a proof
of this assertion in the case $n=2$. Their result appears to be
complete for varieties whose Gauss maps have one-dimensional
fibers. However, their result is incomplete for tangentially
degenerate hypersurfaces whose Gauss maps have fibers of dimension
greater than one. In the recently paper [AG 01b], Akivis and
Goldberg showed that there exist  hypersurfaces which cannot be
built out of cones and osculating varieties.

In the current paper we define and study new types of varieties
with degenerate Gauss maps: varieties with multiple foci and their
particular case, the so-called twisted cones.

\section{Basic equations of a hypersurface of rank $\boldsymbol{r}$
with $\boldsymbol{r}$-multiple focus hyperplanes}

In our recently published paper [AG 02]  in a projective space
$P^N$, we considered varieties $X$ with degenerate Gauss maps of
dimension $n$ and rank $r$ with the following two properties:

\begin{description}
\item[(i)] Their focus hypersurfaces $F_L$ degenerate into
$r$-fold hyperplanes.

\item[(ii)] Their system of second fundamental forms possesses at
least two forms, whose  $\lambda$-equation has $r$ distinct roots.
\end{description}
We have proved that such  varieties $X$ are cones in the space
$P^N$ with a vertex of dimension $l - 1$, where $l = n - r$.

In this paper we also  consider  varieties $X$ with degenerate
Gauss maps of dimension $n$ and rank $r$ with $r$-fold focus
hyperplanes but we  assume that all their second fundamental forms
are proportional, i.e., for each pair of  second fundamental forms
of $X$, their $\lambda$-equation has $r$-multiple eigenvalues.

Since we assume $r \geq 2$,  Segre's theorem (see [AG 93], Theorem
2.2, p. 55) implies that such  varieties are hypersurfaces in a
subspace $P^{n+1}$. We shall prove that such hypersurfaces can
differ from cones.

Consider a  hypersurface $X$ with a degenerate Gauss map of
dimension $n$ and rank $r$ whose focus hypersurfaces $F_L$ are
$r$-fold hyperplanes of dimension $l-1$, where $l = n - r$ is the
dimension of the Monge--Amp\`{e}re foliation on $X$. We associate
a family of moving frames $\{A_u\},  \, u = 0, 1, \ldots , n + 1$,
with $X$ in such a way that the point $A_0$ is a regular point of
a generator $L$, the points $A_a, \, a = 1, \ldots, l$, belong to
the $r$-fold focus hyperplane $F_L$, the points $A_p, \, p = l +
1, \ldots, n$, lie in the tangent hyperplane $T_L (X)$, and the
point $A_{n+1}$ is situated outside of this hyperplane.

The equations of infinitesimal displacement of the moving frame
$\{A_u\}$ are
 \begin{equation}\label{eq:1}
 dA_u = \omega_u^v A_v, \;\;\;\;u, v = 0, 1, \ldots , n + 1,
\end{equation}
where $\omega_u^v$ are 1-forms satisfying the structure equations
of the projective space $P^N$:
 \begin{equation}\label{eq:2}
 d \omega_u^v = \omega_u^w \wedge \omega_w^v, \;\;\;\;u, v, w = 0, 1, \ldots , n +
 1.
\end{equation}

As a result of the specialization of the moving frame mentioned
above, we obtain the following basic equations of the variety $X$:
\begin{equation}\label{eq:3}
 \omega_0^{n+1} = 0, \;\; \omega^{n+1}_a = 0, \;\;\; a = 1,
\ldots, l,
\end{equation}
\begin{equation}\label{eq:4}
 \omega_p^{n+1} = b_{pq} \omega^q, \;\;  \omega_a^p = c_{aq}^p \omega^q,
\end{equation}
and
\begin{equation}\label{eq:5}
b_{sq} c^s_{ap} = b_{sp} c^s_{aq};
\end{equation}
here $\omega^q := \omega_0^q$ are the basis forms of the variety
$X$, and $B = (b_{pq})$ is a nondegenerate symmetric $(r \times
r)-$ matrix (see [AG 93], Section 4.1).

Denote by $C_a$  the $(r \times r)$-matrix occurring in equations
(4):
$$
C_a = ( c^p_{aq}).
$$
If we  use the identity matrix $C_0 = (\delta^p_q)$ and the index
$i = 0, 1, \ldots, l$ (i.e., $\{i\} = \{0, a\}$), then equations
$b_{pq} = b_{qp}$ and (5) can be combined and written as follows:
\begin{equation}\label{eq:6}
(B C_i)^T = (B C_i),
\end{equation}
i.e., the matrices
$$
H_i = B C_i = (b_{qs} c^s_{ip})
$$
are symmetric.

Since the points $A_a, \, a = 1, \ldots, l$, belong to
the $r$-fold focus $(l-1)$-plane $F_L$,  equation of $F_L$ is
$$
(x^0)^r = 0.
$$
But in the general case the focus hypersurface $F_L$ of the
generator $L$ is determined by the equation $$ \det \, (\delta_q^p
x^0 + c^p_{aq} x^a) = 0 $$ (see equation (4.19) on p. 117
of the book [AG 93]).
Thus, we have
$$
   \det (\delta_q^p x^0 + c^p_{aq} x^a) = (x^0)^r
$$
It follows that each of the matrices $C_a$ has an $r$-multiple
eigenvalue $0$, and as a result, each of these matrices is
nilpotent. We assume that each of the matrices $C_a$ has the form
\begin{equation}\label{eq:7}
C_a = (c_{aq}^p), \; \text{where} \; c^p_{aq} = 0 \; \text{for} \;
p \geq q.
\end{equation}
Thus,  $\text{rank} \; C_a \leq r - 1$. It follows that all
matrices $C_a$ are nilpotent. Denote by $r_1$ the maximal rank of
matrices from the bundle $C = x^a C_a, \, r_1 \leq r - 1$.

It is obvious that this form is sufficient for all $F_L$ to be
$r$-fold hyperplanes. Wu and Zheng [WZ 02] (see also Piontkowski
[P 01, 02]) proved this for the ranks $r = 2, 3, 4$ and different
values of the maximum rank $r_1$ of matrices of the bundle $x^a
C_a$. However, Wu and Zheng in [WZ 02] gave also a counterexample
which proves that the form (7) is not necessary for all $F_L$ to
be $r$-fold hyperplanes.

A single second fundamental form of $X$ at its regular point $x =
A_0$ can be written as
$$
\Phi_0 = b_{pq} \omega^p \omega^q.
$$
This form is of rank $r$. At the singular points $A_a$ belonging
to an $r$-multiple focus hyperplane $F_L$, the second  fundamental
form of the hypersurface $X$ has the form
\begin{equation}\label{eq:8}
\Phi_a = b_{ps} c^s_{aq} \omega^p \omega^q,
\end{equation}
where $(b_{ps} c^s_{aq})$ is a symmetric matrix. The maximal rank
of matrices from the bundle $\Phi = x^a \Phi_a$ is also equal to
$r_1 < r - 1$.

\section{Hypersurfaces with degenerate Gauss maps of rank
$\boldsymbol{r}$ with a one-dimensional Monge--Amp\`{e}re
foliation and $\boldsymbol{r}$-multiple foci} Let $A_0 A_1$ be a
leaf of the Monge--Amp\'{e}re foliation, let $A_0$ be a regular
point of this foliation, and let $A_1$ be its $r$-multiple focus.
Then in equations (6), we have $a, b = 1; \, p, q = 2, \ldots, n$,
and these equations become
\begin{equation}\label{eq:9}
 \omega_p^{n+1} = b_{pq} \omega^q, \;\;  \omega_1^p = c_q^p
 \omega^q.
\end{equation}

By our assumption (7), the matrix $C = (c_q^p)$ has the form
\begin{equation}\label{eq:10}
\renewcommand{\arraystretch}{1.3}
 C =
\left(
\begin{array}{cccc}
0 & c_3^2 &\ldots & c_n^2\\
\hdotsfor4  \\
0 & 0  & \ldots      & c_n^{n-1}\\
0 & 0  &\ldots       & 0
\end{array}
\right),
\renewcommand{\arraystretch}{1}
\end{equation}
where the coefficients $c_{p+1}^p \neq 0$. As to the matrix $B =
(b_{pq})$, by the relation
\begin{equation}\label{eq:11}
BC = CB
\end{equation}
(cf. (6)), this matrix has the form
\begin{equation}\label{eq:12}
 \renewcommand{\arraystretch}{1.3}
 B =
\left(
\begin{array}{cccc}
0 & \ldots & 0 & b_{2,n}\\
0 & \ldots & b_{3, n-1} & b_{3,n}\\
\hdotsfor4  \\
b_{n,2} & \ldots & b_{n,n-1} & b_{nn}
\end{array}
\right),
\renewcommand{\arraystretch}{1}
\end{equation}
and $\text{rank} \; C = n - 2, \, \text{rank} \; B = n - 1$. In
addition, by (11), the entries of  the matrices $B$ and $C$ are
connected by certain bilinear relations implied by (11).

By (9), (10), and (12), on the hypersurface $X$, we have the
equation
\begin{equation}\label{eq:13}
\omega_1^n = 0.
\end{equation}
Since on the hypersurface $X$ also the equations (3) hold, the
differentials of the points $A_0$ and $A_1$ take the form
\begin{equation}\label{eq:14}
  \renewcommand{\arraystretch}{1.3}
 \begin{array}{ll}
dA_0 = \omega_0^0 A_0 + \omega_0^1 A_1 + \omega_0^2 A_2 +
     \ldots + \omega_0^{n-1} A_{n-1} + \omega_0^n A_n,\\
dA_1 = \omega_1^0 A_0 + \omega_1^1 A_1 + \omega_1^2 A_2 +
    \ldots + \omega_1^{n-1} A_{n-1}.
\end{array}
\renewcommand{\arraystretch}{1}
\end{equation}
In equations (14), the forms $\omega_1^2, \omega_1^3, \ldots ,
\omega_1^{n-1}$ are linearly independent, and by (9) and (10),
they are expressed in terms of the basis forms $\omega^3, \ldots,
\omega^n$ only. The following cases can occur:

\begin{description}
\item[1)] The 1-form $\omega_1^0$ is independent of the forms
$\omega^3, \ldots, \omega^n$, and hence also of the forms
$\omega_1^2, \ldots, \omega_1^{n-1}$. In this case, the
$r$-multiple focus $A_1$ of the rectilinear generator $L$
describes a focus variety $G$ of dimension $r = n - 1$. The
variety $G$ is of codimension two in the space $P^{n+1}$ in which
the hypersurface $X$ is embedded. The tangent subspace $T_{A_1}
(G)$ is defined by the points $A_1, A_0, A_2, \ldots , A_{n-1}$.
At the point $A_1$, the variety $G$ has two independent second
fundamental forms. We can determine these two forms by finding the
second differential of the point $A_1$:
 $$
d^2 A_1 \equiv \omega_1^p \omega_p^n A_n +  \omega_1^p
\omega_p^{n+1} A_{n+1} \pmod{T_{A_1}(G)}.
 $$
 Thus, we have
 $$
\Phi_1^n = \omega_1^p \omega_p^n, \;\; \Phi_1^{n+1} = \omega_1^p
\omega_p^{n+1}.
 $$
The second of these forms coincides with the second fundamental
form $\Phi_1$ of the hypersurface $X$ at the point $A_1$. By (10),
if $\omega^3 = \ldots = \omega^n = 0$, the 1-forms $\omega_1^p =
0$. Hence the quadratic forms $\Phi_1^{n}$ and $\Phi_1^{n+1}$
vanish on the focal variety $G$. Thus the direction $A_1 \wedge
A_0$ is an asymptotic direction on $G$.

\item[2)] The 1-form $\omega_1^0$ is a linear combination  of the
forms  $\omega_1^2, \ldots, \omega_1^{n-1}$, and hence also of the
forms  $\omega^3, \ldots, \omega^n$. In this case, the focus $A_1$
of the rectilinear generator $L$ describes a focus  variety $G$ of
dimension $n - 2$, and its  tangent subspace $T_{A_1} (G)$ is a
hyperplane in the space $A_0 \wedge A_1 \wedge A_2 \wedge \ldots
\wedge A_{n-1}$. For $\omega_1^2 = \ldots = \omega_1^{n-1} = 0$,
the point $A_1$ is fixed, and the straight line $L = A_1 A_0$
describes a two-dimensional cone with vertex $A_1$. This cone is
called the \emph{fiber cone}. The hypersurface $X$ foliates into
an $(n-2)$-parameter family of such fiber cones. It is called a
\emph{twisted cone with rectilinear generators}.

Later on, in Section 3, for $n = 3$ we will prove that a fiber
cone is a pencil of straight lines. Most likely this is true for
any $n$.

\item[3)] Suppose that an $(n-2)$-dimensional focus  variety $G$
of the hypersurface $X$ belongs to a hyperplane $P^n$ of the space
$P^{n+1}$. We can take this hyperplane as the hyperplane at
infinity $P^n_\infty$ of  the space $P^{n+1}$. As a result, the
space $P^{n+1}$ becomes an affine space $A^{n+1}$. In this case,
the hypersurface $X$ becomes also a \emph{twisted cylinder} in
$A^{n+1}$, which foliates into an $(n-2)$-parameter family of
two-dimensional cylinders with rectilinear generators. The
hypersurface $X$ with a degenerate Gauss map is not a cylinder in
$A^{n+1}$ and does not have singularities in this space. Thus,
this hypersurface is an affinely complete hypersurface in
$A^{n+1}$, which is not a cylinder. An example of such a
hypersurface in the space $A^{n+1}$ was considered by Sacksteder
and Bourgain (see [S 60], [W 95], [I 98, 99a, b], and [AG 01a]).
  \end{description}

Note also that  hypersurfaces  with  degenerate Gauss maps in the
space $P^{n+1}$ considered in this section are \emph{lightlike
hypersurfaces} which were studied in detail in the papers [AG 98a,
b; 99a, b] by Akivis and Goldberg.

\section{Hypersurfaces
with  Degenerate Gauss Maps with Double Foci on Their Rectilinear
\newline Generators  in the  Space $\boldsymbol{P^4}$}

As an example, we consider hypersurfaces $X$ with degenerate Gauss
maps of rank $r = 2$ in the space $P^4$ that have a double focus
$F$ on each rectilinear generator $L$. On a rectilinear generator
$L = A_0 \wedge A_1$ of $X$, there is a single (double) focus
$F_L$. With respect to a first-order frame, the basic equations of
$X$ are
\begin{equation}\label{eq:15}
\omega_0^4 = 0, \;\; \omega_1^4 = 0.
\end{equation}
The basis forms of $X$ are $\omega_0^2$ and $\omega_0^3$. By
(4) and (12), with respect to a second-order frame, we have
the following equations
\begin{equation}\label{eq:16}
 \renewcommand{\arraystretch}{1.3}
 \left\{
 \begin{array}{lll}
 \omega_2^4 =  &\!\!\!\!b_{23} \omega^3_0, & \omega_1^2 = c_3^2 \omega^3_0, \\
\omega_3^4 = b_{32} \omega^2_0 + &\!\!\!\!b_{33} \omega^3_0, &
\omega_1^3 = 0,
\end{array}
\right.
\renewcommand{\arraystretch}{1}
\end{equation}
where $b_{23} = b_{32} \neq 0$ and $c_3^2 \neq 0$. As a result,
the matrices $B$ and $C$ take the form
$$
B = \left( \begin{array}{cc}
0 & b_{23} \\
b_{23} & b_{33}
\end{array}
\right), \;\; C = \left( \begin{array}{cc}
0 & c_3^2 \\
0 & 0
\end{array}
\right).
$$

The differentials of the points $A_0$ and $A_1$ are
$$
  \renewcommand{\arraystretch}{1.3}
 \begin{array}{ll}
dA_0 = \omega_0^0 A_0 + \omega_0^1 A_1 + \omega_0^2 A_2
      + \omega_0^{3} A_{3},\\
dA_1 = \omega_1^0 A_0 + \omega_1^1 A_1 + \omega_1^2 A_2.
    \end{array}
\renewcommand{\arraystretch}{1}
$$
(cf. (14)). The point $A_1 = F_L$ is a single focus of a
rectilinear generator $L$.

Exterior differentiation of equations (16) gives the following
exterior quadratic equations:
\begin{equation}\label{eq:17}
- 2 b_{23} \omega_2^3 \wedge \omega_0^2 + \Delta b_{23}  \wedge
\omega_0^3  = 0,
\end{equation}
\begin{equation}\label{eq:18}
\Delta b_{23} \wedge \omega_0^2 + \Delta b_{33}  \wedge \omega_0^3
= 0,
\end{equation}
\begin{equation}\label{eq:19}
- (\omega_1^0 + c_3^2 \omega_2^3) \wedge \omega_0^2 + \Delta
c_{3}^2 \wedge \omega_0^3  = 0,
\end{equation}
\begin{equation}\label{eq:20}
(\omega_1^0 - c_3^2 \omega_2^3)  \wedge \omega_0^3  = 0,
\end{equation}
where
$$
  \renewcommand{\arraystretch}{1.3}
 \begin{array}{ll}
\Delta b_{23} = db_{23} + b_{23} (\omega_0^0  - \omega_2^2 -
\omega_3^3
 + \omega_4^{4}) - b_{33} \omega_2^3, \\
 \Delta b_{33} = db_{33} + b_{33} (\omega_0^0 - 2\omega_3^3 + \omega_4^{4})
 + b_{32} c_3^2  \omega_0^1 - b_{32} \omega_3^2,\\
 \Delta c_{3}^2 = dc_{3}^2 + c_{3}^2 (\omega_0^0 - \omega_1^1 + \omega_2^2
 - \omega_3^{3}).
    \end{array}
\renewcommand{\arraystretch}{1}
$$
From equations (17) and (20), it follows that the forms
$\omega_2^3$ and $\omega_1^0$ are linear combinations of the basis
forms $\omega_0^2$ and $\omega_0^3$. Three cases are possible:

\begin{description}

\item[1)] $\omega_1^0 \wedge \omega_0^3 \neq 0$. Since by (16),
this implies that $\omega_1^0 \wedge \omega_1^2 \neq 0$, it
follows that the focus $A_1$ describes a two-dimensional focal
surface $G^2$. The tangent plane to $G^2$ at the point $A_1$ is
$T_{A_1} (G) = A_1 \wedge A_0 \wedge A_2$, and the straight line
$L = A_0 A_1$ is tangent to $G^2$ at $A_1$.

\item[2)]
\begin{equation}\label{eq:21}
\omega_1^0 \wedge \omega_0^3  = 0.
\end{equation}
In this case, the point $A_1$ describes a focal line $G^1$,  and
the straight line $L = A_0 A_1$ intersects this line $G^1$ at the
point $A_1$. The hypersurface $X$ foliates into a one-parameter
family of two-dimensional cones and is a twisted cone.

\item[3)] The osculating hyperplane of the curve $G^1$ described
in 2) is fixed.
\end{description}

We consider these three cases in detail.
\begin{description}

\item[1)] We prove an existence theorem for this case applying the
Cartan test (see, for example, [BCGGG 91]).

\begin{theorem} Hypersurfaces $X$ of rank two
in the space $P^4$, for which the single focus of a rectilinear
generator $L$ describes a two-dimensional surface, exist, and the
general solution of the system defining such hypersurfaces depends
on one function of two variables. The direction $A_1 A_0$ is an
asymptotic direction on the surfaces $G^2$, and the hypersurface
$X$ is formed by the asymptotic tangents to the surfaces $G^2$.
\end{theorem}

\begin{proof} On a hypersurface in question, the inequality
$\omega^1_0 \wedge \omega^3_0 \neq 0$ holds as well as the
exterior quadratic equations (17)--(20). The latter equations
contain five forms $\omega_2^3,\, \Delta b_{23},\, \Delta b_{33},
\,\omega_1^0$, and $\Delta c_{3}^2$ that are different from the
basis forms $\omega_0^2$ and $\omega_0^3$. So, we have $q = 5$.

The character $s_1$ of the system under investigation is equal to
the number of independent exterior quadratic equations
(17)--(20). Thus, we have $s_1 = 4$. As a result, the second
character of the system is $s_2 = q - s_1 = 1$. Therefore, the
Cartan number $Q = s_1 + 2 s_2 = 6$.

We calculate now the number $S$ of parameters on which the most
general integral element of the system under investigation depends
(i.e., the dimension $S$ of the space of integral elements over a
point). Applying Cartan's lemma to equations (17) and (18), we
find that
\begin{equation}\label{eq:22}
 \renewcommand{\arraystretch}{1.3}
 \left\{
 \begin{array}{rll}
- 2b_{23} \omega_2^3 &\!\!\!\!= b_{222} \omega^2 + b_{223}\omega^3, \\
\Delta b_{23} &\!\!\!\! = b_{232} \omega^2 + b_{233}\omega^3, \\
\Delta b_{33} &\!\!\!\!= b_{332} \omega^2 + b_{333}\omega^3_0.
\end{array}
\right.
\renewcommand{\arraystretch}{1}
\end{equation}
Since the coefficients of the basis forms in the right-hand sides
of (22) are symmetric with respect to the lower indices, the
number of independent among these coefficients is $S_1 = 4$.

Equation (20) implies that
\begin{equation}\label{eq:23}
\omega_1^0 = c_3^2 \omega_2^3  + \lambda \, \omega^3_0.
\end{equation}
We substitute this expression into equation (19). As a result,
we obtain
\begin{equation}\label{eq:24}
 - 2(c_3^2 \omega_2^3  + \lambda \,\omega_0^3)  \wedge \omega_0^2
 + \Delta c_{3}^2 \wedge \omega_0^3  =0.
\end{equation}
It follows from (24) that the 1-form $\Delta c_{3}^2$ is a
linear combination of the basis forms. We write this expression in
the form
\begin{equation}\label{eq:25}
 \Delta c_{3}^2 = \mu \, \omega_0^2 + \nu \, \omega_0^3.
\end{equation}
We can find the form $\omega_2^3$ from the first equation of
(22). Substituting this expression and (25) into equation
(24), we find that
$$
\Biggl(\displaystyle \frac{c_3^2 b_{223}}{b_{23}} - \lambda
\Biggr)  \omega_0^3 \wedge \omega_0^2 + \mu \, \omega_0^2 \wedge
\omega_0^3= 0.
$$
This implies that
$$
\mu = \displaystyle \frac{c_3^2 b_{223}}{b_{23}} - \lambda.
$$
Thus, there are only two independent coefficients in
decompositions (23) and (25), $S_2 = 2$. As a result, we have $S =
S_1 + S_2 = 6$, and $S = Q$. Thus, by Cartan's test, the system
under investigation is in involution, and its general solution
depends on one function of two variables.

Next, we find the second fundamental forms of the two-dimensional
focal surface $G^2$ of the hypersurface $X$ with a degenerate
Gauss map. To this end, we compute
$$
d^2 A_1 \equiv (\omega_1^0 \omega_0^3 + \omega_1^2 \omega_2^3) A_3
+ \omega_1^2 \omega_2^4 A_4 \pmod{T_{A_1} (G^2)}.
$$
Thus, the second fundamental forms of $G^2$ are
$$
\Phi_1^3 = \omega_1^0 \omega_0^3 + \omega_1^2 \omega_2^3, \;\;
\Phi_1^4 = \omega_1^2 \omega_2^4.
$$
The direction $A_1 A_0$ is defined on $G_2$ by the equation
$\omega_1^2 = 0$. By (16), this equation is equivalent to the
equation $\omega^3_0 = 0$. Thus, in this direction the second
fundamental forms $\Phi_1^3$ and $\Phi_1^4$ vanish:
$$
\Phi_1^3 \equiv 0 \pmod{\omega^3_0}, \;\; \Phi_1^4 \equiv 0
\pmod{\omega^3_0},
$$
and the direction $A_1 A_0$ is an asymptotic direction on the
focal surface $G^2$. \hfill
\end{proof}

\item[2)] We prove the following existence theorem for the twisted
cones.
\begin{theorem} If condition $(21)$ is satisfied, then the double focus $A_1$
of the generator $A_0 \wedge A_1$ of the variety $X$ describe the
focal curve, and $X$ is a twistor cone. In the space $P^4$, the
twisted cones exist, and the general solution of the system
defining such cones depends on five functions of one variable.
\end{theorem}

\begin{proof}
In this case, the point $A_1$ describes the focal line $G^1$. Thus
we must enlarge the system of equations (16) by the equation
\begin{equation}\label{eq:26}
 \omega_1^0 = a \,\omega_0^3.
\end{equation}
Equation (26) is equivalent to equation (21). The 1-form
$\omega_0^3$ is a basis form on the focal line $G^1$. By (26),
equation (20) takes the form
$$
 \omega_2^3 \wedge \omega_0^3 = 0.
$$
It follows that
\begin{equation}\label{eq:27}
 \omega_2^3 = b\, \omega_0^3.
\end{equation}
Now equations (27) and (19) become
\begin{equation}\label{eq:28}
(\Delta b_{23} + 2 b_{23}\, b \,\omega_0^2)  \wedge \omega_0^3  =
0,
\end{equation}
\begin{equation}\label{eq:29}
(\Delta c_{3}^2 + (a + b \, c_3^2) \,\omega^2) \wedge \omega_0^3
= 0.
\end{equation}
Equation (28) remains the same.

Taking exterior derivatives of equations (26) and (28), we
obtain the exterior quadratic equations
\begin{equation}\label{eq:30}
(da + a (2  \omega_0^0 - \omega_1^1 - \omega_3^3) + c_3^2 \,
\omega_2^0 + ab \,\omega_0^2)  \wedge \omega_0^3 = 0,
\end{equation}
\begin{equation}\label{eq:31}
(db + b (\omega_0^0 - \omega_2^2)+ b_{23}\, \omega_4^3 + b
\,\omega_0^2) \wedge \omega_0^3 = 0.
\end{equation}
Now the system of exterior quadratic equations consists of
equations (18), (28)--(31). Thus, we have $s_1 = 5$. In addition
to the basis forms $\omega_0^2$ and $\omega_0^3$, these exterior
equations contain the forms $\Delta b_{23}, \, \Delta b_{33},
\Delta c_{3}^2, \, \Delta a$, and $\Delta b$, where
\begin{equation}\label{eq:32}
\Delta a = da + a (2  \omega_0^0 - \omega_1^1 - \omega_3^3) +
c_3^2 \omega_2^0
\end{equation}
and
$$
\Delta b = db + b (\omega_0^0 - \omega_2^2)+ b_{23} \omega_4^3.
$$
The number of these forms is $q = 5$. Thus, $s_2 = q - s_1 = 0$,
and the Cartan number $Q = s_1 = 5$. If we find the forms $\Delta
b_{23}, \, \Delta b_{33},\, \Delta c_{3}^2, \,\Delta a$, and
$\Delta b$ from  the system of equations (18), (28)--(31), we see
that the most general integral element of the system under
investigation (i.e., the dimension $S$ of the space of integral
elements over a point) depends on $S = 5$ parameters. Thus, $S =
Q$, the system under investigation is in involution, and its
general solution depends on five functions of one variable. \hfill
\end{proof}

Consider the focal curve $G^1$ of the twisted cone $X^3 \subset
P^4$ described by the point $A_1$. We have
$$
dA_1 = \omega_1^1 A_1 + (c_3^2 A_2 + a A_0)\, \omega_0^3.
$$
The point $\widetilde{A}_2 = c_3^2 A_2 + a A_0$ along with the
point $A_1$ define a tangent line to $G^1$. We can specialize our
moving frame by locating its vertex $A_2$ at $\widetilde{A}_2 $
and by normalizing the frame by means of the condition $c_3^2 =
1$. Then we obtain
$$
dA_1 = \omega_1^1 A_1 + \omega_0^3 A_2.
$$
In addition, the conditions
$$
a = 0, \;\; c_3^2 = 1
$$
are satisfied. These conditions and equations (16), (26),
(29), and (32) imply that
\begin{equation}\label{eq:33}
 \omega_1^2 = \omega_0^3, \;\; \omega_1^0 = 0,
\end{equation}
\begin{equation}\label{eq:34}
\Delta c_{3}^2 = \omega_0^0 - \omega_1^1 + \omega_2^2 -\omega_3^3,
\end{equation}
\begin{equation}\label{eq:35}
\Delta a = \omega_2^0.
\end{equation}

After this specialization, the straight line $A_1 \wedge A_2$
becomes the tangent to the focal line $G^1$.

Now equations (29) and (30) take the form
$$
\renewcommand{\arraystretch}{1.3}
  \begin{array}{ll}
 (\omega_0^0 - \omega_1^1 + \omega_2^2 -\omega_3^3 + b \,\omega_0^2)
 \wedge \omega_0^3 = 0, \\
 \omega_2^0  \wedge \omega_0^3 = 0.
\end{array}
$$
It follows from the last equation that
\begin{equation}\label{eq:36}
 \omega_2^0 = c \, \omega_0^3.
\end{equation}

Note also that equation (31) shows that since $b_{23} \neq 0$, the
quantity $b$ can be reduced to 0 by means of the form
$\omega_4^3$. As a result, equation (27) takes the form
\begin{equation}\label{eq:37}
 \omega_2^3 = 0,
\end{equation}
and equation (31) becomes
\begin{equation}\label{eq:38}
 \omega_4^3 \wedge \omega_0^3 = 0.
\end{equation}
It follows from (38) that
\begin{equation}\label{eq:39}
 \omega_4^3 = f \,\omega_0^3.
\end{equation}

Differentiating the point $A_2$ and applying (16), (36),  and
(37), we obtain
$$
dA_2 = \omega_2^2 A_2 + \omega_2^1 A_1 + (c A_0 + b_{23} A_4)\,
\omega_0^3.
$$
The 2-plane $\alpha = A_1 \wedge A_2 \wedge (c A_0 + b_{23} A_4)$
is the \emph{osculating plane} of the line $G^1$ at the point
$A_1$.

We place the point $A_4$ of our moving frame into the plane
$\alpha$ and make a normalization $b_{23} = 1$. As a result, we
have $c = 0$ and
\begin{equation}\label{eq:40}
 \omega_2^0 = 0, \;\;\omega_2^4 = \omega_0^3.
\end{equation}
Now, the plane $\alpha$ is defined as
 $\alpha = A_1 \wedge A_2 \wedge A_4$,
and the differential of $A_2$ becomes
$$
dA_2 = \omega_2^2 A_2 + \omega_2^1 \, A_1 +  \omega_0^3 \, A_4.
$$

Taking the exterior derivative of the first equation of (40), we
obtain
$$
\omega_4^0 \wedge \omega_0^3 = 0,
$$
and this implies that
\begin{equation}\label{eq:41}
 \omega_4^0 = g\, \omega_0^3.
\end{equation}

By means of equations (37) and (41), we find that
\begin{equation}\label{eq:42}
d A_4 = \omega_4^4 A_4 +  \omega_4^1 A_1 +
 \omega_4^2 A_2  + (f A_3 + g A_0) \, \omega_0^3.
\end{equation}
Equation (42) means that the 3-plane
$$
\beta = A_1 \wedge A_2 \wedge A_4 \wedge (f A_3 + g A_0)
$$
is the \emph{osculating hyperplane} of the focal line $G^1$.

Taking exterior derivatives of equations (39) and (41), we
find the following exterior quadratic equations:
\begin{equation}\label{eq:43}
(df + f (\omega_0^0 - \omega_4^4)) \wedge \omega_0^3 = 0,
\end{equation}
and
\begin{equation}\label{eq:44}
(dg + g (2\omega_0^0 - \omega_3^3 - \omega_4^4) - f \omega_3^0)
\wedge \omega_0^3 = 0.
\end{equation}
As we did earlier, we can prove that by means of the secondary
forms $\omega_0^0 - \omega_4^4$ and $\omega_3^0$, we can
specialize our moving frames in such a way, that
$$
f = 1, \;\; g = 0.
$$
As a result, equations (39) and (41) become
\begin{equation}\label{eq:45}
\omega_4^3 = \omega_0^3, \;\; \omega_4^0 = 0,
\end{equation}
and the osculating hyperplane $\beta$ of $G^1$ becomes
$$
\beta = A_1 \wedge A_2 \wedge A_4 \wedge A_3.
$$

Substituting the values $f = 1$ and $g = 0$ into equations (43)
and (44), we obtain
\begin{equation}\label{eq:46}
(\omega_0^0 - \omega_4^4) \wedge \omega_0^3 = 0,
\end{equation}
and
\begin{equation}\label{eq:47}
 \omega_3^0 \wedge \omega_0^3 = 0.
\end{equation}
Note that equations (46) and (47) could be also obtained by
exterior differentiation of equations (45).

After the above specialization, we obtain the following system of
equations defining the twisted cones $X$ in the space $P^4$:
\begin{equation}\label{eq:48}
 \renewcommand{\arraystretch}{1.3}
 \left\{
 \begin{array}{ll}
\omega_2^4 = \omega_0^3, & \omega_3^4 = \omega^2, \\
\omega_1^2 = \omega_0^3, & \omega_1^3 = 0, \\
\omega_1^0 = 0,        & \omega_2^3 = 0, \\
\omega_2^0 = 0,         & \omega_2^4 = \omega^3, \\
\omega_4^3 = \omega_0^3,  & \omega_4^0 = 0.
\end{array}
\right.
\renewcommand{\arraystretch}{1}
\end{equation}
Note that in addition to all specializations made earlier, in
equations (48), we also made a specialization $b_{33} = 0$ which
can be achieved by means of the secondary form $\omega_0^1 -
\omega_3^2$ (see the third equation in (22)).

Taking exterior derivatives of equations (48),
  we find the following exterior quadratic equations:
\begin{equation}\label{eq:49}
 \renewcommand{\arraystretch}{1.3}
 \left\{
 \begin{array}{ll}
(\omega_0^0 - \omega_2^2 - \omega_3^3 + \omega_4^4) \wedge \omega_0^3 = 0,  \\
(\omega_0^0 - \omega_2^2 - \omega_3^3 + \omega_4^4) \wedge
\omega_0^2
+ (\omega_0^1 - \omega_3^2)  \wedge \omega_0^3 = 0,  \\
(\omega_0^0 - \omega_1^1 + \omega_2^2 - \omega_3^3) \wedge \omega_0^3 = 0,  \\
(\omega_0^0 - \omega_4^4) \wedge \omega_0^3 = 0,  \\
\omega_3^0 \wedge \omega_0^3 = 0.
\end{array}
\right.
\renewcommand{\arraystretch}{1}
\end{equation}
The exterior differentiation of the remaining five equations of
system (48) leads to identities.

The system of equations (49) is equivalent to the system of
equations (18), (28)--(31) from which it is obtained as a result
of specializations. For the system of equations (49), as well as
for the original system  of equations (18), (28)--(31), we have $q
= 5, \, s_1 = 5, \, s_2 = 0, \, Q = N = 5$. \emph{The system is in
involution, and its general solution exists and depends on five
functions of one variable.}

Further we investigate the structure of the fiber cones of a
twisted cone $X \subset P^4$. The fiber cones on $X$ are defined
by the system
\begin{equation}\label{eq:50}
  \omega_0^3 = 0.
\end{equation}
By (50) and (48), we have
\begin{equation}\label{eq:51}
d A_0 =  \omega_0^0 A_0 +  \omega_0^1 A_1 + \omega_0^2 A_2.
\end{equation}
It follows that the plane $A_0 \wedge A_1 \wedge  A_2$ is tangent
to the fiber cone $C$ along its generator $L = A_0 \wedge A_1$.
 By (50) and (48), the differential of the point $A_2$ is
\begin{equation}\label{eq:52}
d A_2 =  \omega_2^1 A_1 +  \omega_2^2 A_2,
\end{equation}
and by (50), we also have
\begin{equation}\label{eq:53}
d A_1 =  \omega_1^1 A_1.
\end{equation}
Equations (51), (52), and (53) prove that the tangent plane
$\gamma = A_0 \wedge A_1 \wedge A_2$ to the fiber cone $C$ is
fixed when the generator $L = A_0 \wedge A_1$ moves along $C$. It
follows that a  fiber cone $C$ is simply a pencil of straight
lines with center at $A_1$ located at the plane $\gamma$.

Thus, we have proved the following theorem.

\begin{theorem}
A twisted cone $X$ in the space $P^4$ foliates into a
one-parameter family of pencils of  straight lines whose centers
are located on its focal line $G^1$ and whose planes are tangent
to $G^1$.
\end{theorem}

Namely such a picture can be seen in the example of
Sacksteder--Bourgain (see [AG 01a]). However, here we proved this
theorem for the general case.

Now we prove the converse statement: A general smooth
one-parameter family of two-dimensional planes $\gamma (t)$ in the
space $P^4$ forms a three-dimensional twisted cone $X$. In fact,
such a family envelopes a curve $G^1$, whose point $A$ is the
common point of the planes $\gamma (t)$ and  $\gamma (t + dt)$,
i.e., $A (t) =  \gamma (t) \cap \gamma (t+dt)$. The point $A (t)$
and the plane  $\gamma (t)$ define a pencil $(A,  \gamma) (t)$ of
straight lines with center $A (t)$ and plane  $\gamma (t)$. The
set of these pencils forms a  three-dimensi\-onal ruled surface
$X$ with rectilinear generators $L$ belonging to the pencils $(A,
\gamma) (t)$. Moreover, the tangent space $T (X)$ is constant
along a rectilinear generator $L$. Hence the rank of the
 variety $X$ equals two.

Since the dimension of the Grassmannian $G (2, 4)$ of
two-dimensional planes in the space $P^4$ is equal to six (see [AG
93], Section 1.4, p. 297), one-parameter family of such planes
depends on five functions of one variable. This number coincides
with the arbitrariness of existence of twisted cones in $P^4$ that
we computed earlier by investigating a system defining a twisted
cone (see Theorem 2).

\item[3)] Next we find under what condition a twisted cone becomes
a twisted cylinder. This condition is equivalent to a condition
under which the osculating hyperplane $\beta$ of the focal curve
$G^1$ is fixed, when the point $A_1$ moves along $G^1$. Since
$\beta = A_1 \wedge A_2 \wedge A_3 \wedge A_4$ and
$$
d A_3 =  \omega_3^0 A_0 +  \omega_3^1 A_1 + \omega_3^2 A_2
 +  \omega_3^3 A_3 + \omega_3^4 A_4,
$$
the condition in question has the form
\begin{equation}\label{eq:54}
 \omega_3^0 = 0.
\end{equation}
If we take the fixed osculating hyperplane $\beta$ of $G^1$ as the
hyperplane at infinity $H_\infty$ of the space $P^4$, then $P^4$
becomes an affine space $A^4$. Then the hypersurface $X$ becomes a
twisted cylinder $\widetilde{X}$, which by Theorem 3,
\emph{foliates into a one-parameter family of planar pencils of
parallel straight lines.} The hypersurface $X$ does not have
singularities in the space $A^4$ and is a \emph{complete smooth
noncylindrical hypersurface.}

It is easy to prove the existence of twisted cylinders  in the
space $A^4$.

\begin{theorem}
Twisted cylinders  in the space $A^4$ exist, and the general
solution of the system defining such cylinders depends on four
functions of one variable.
\end{theorem}

\begin{proof}
In fact, a twisted cylinder in $A^4$ is defined by the system of
equations (48) and (54). By (54), the last equation of (49)
becomes an identity. Exterior differentiation of (54) leads to an
identity too. Thus, in the system of exterior quadratic equations
(51), only four equations are independent. Thus, $s_1 = 4$, and
equations (49) contain only four 1-forms that are different from
the basis forms. Hence $q = 4$. Therefore, $s_2 = q - s_1 = 0, \,
Q = s_1 + 2 s_2 = 4$. Equations (49) imply also that $S = 4$.
Since $Q = S$, the system is in involution, and its general
solution depends on four functions of one variable. \hfill
\end{proof}

In conclusion, we indicate a construction defining the general
twisted cylinders in the space $A^4$. Let $P^3$ be an arbitrary
hyperplane in the projective space $P^4$, and let $G^1$ be an
arbitrary curve in $P^3$. Consider a family of planes $\gamma
(t)$, that are tangent to the curve $G^1$ but do not belong to
$P^3$, such that two infinitesimally close planes $\gamma (t)$ and
$\gamma (t+dt)$ of this family do not belong to a
three-dimensional subspace of the space $P^4$. Then these two
planes have only one common point $A (t) =  \gamma (t) \cap \gamma
(t+dt)$ belonging to $G^1$, and the planes $\gamma (t)$ form a
twisted cone in the space $P^4$.  If we take the hyperplane $P^3$
as the hyperplane at infinity of $P^4$, then the space $P^4$
becomes an affine space $A^4$, and a twisted cone formed by the
planes $\gamma (t)$ becomes a twisted cylinder in $A^4$. Such a
construction was considered by Akivis in his paper [A 87].
\end{description}

\noindent {\em Authors' addresses}:\\

\noindent
\begin{tabular}{ll}
M.~A. Akivis &V.~V. Goldberg\\
Department of Mathematics &Department of Mathematical Sciences\\
Jerusalem College of Technology---Mahon Lev &  New
Jersey Institute of Technology \\
Havaad Haleumi St., P. O. B. 16031 & University Heights \\
Jerusalem 91160, Israel &  Newark, N.J. 07102, U.S.A. \\
 & \\
 E-mail address: akivis@mail.jct.ac.il & E-mail address:
 vlgold@m.njit.edu
 \end{tabular}
\end{document}